\documentclass[12pt]{article}
\usepackage{graphicx} % Required for inserting images
\usepackage[utf8]{inputenc}
\usepackage{amsfonts}
\usepackage{amsmath}
\usepackage{amssymb}
\usepackage{setspace}
\usepackage{CJK}
\usepackage{moresize}
\usepackage[utf8]{inputenc}
\usepackage[english]{babel}
\usepackage{xcolor}

    \newcommand{\R}{\mathbb{R}}

\title{A fundamental theorem of algebra for sums of exponentials\footnote{2022 Mathematics Subject Classification: 26Cxx, 37M10, 42A70}}
\author{Pierce Ellingson{\footnote{Department of Radiology, Division of Medical Physics, University of Minnesota, Minneapolis, MN 55455}} \ and Farhad Jafari$^{\dagger,}${\footnote{Emeritus Professor, Department of Mathematics and Statistics, University of Wyoming, Laramie, WY 82070}} }

\date{}

\begin{document}

\maketitle

\section{Statement of the Problem.}

A basic consequence of the fundamental theorem of algebra is that a real (or complex) polynomial of degree $ n $ is uniquely determined from its values at $ n + 1 $ points (counting multiplicity) on the real line (or in the complex plane). Hence a polynomial of degree $  n $ may be sampled at an arbitrary set of $ n+1 $ distinct points to be fully determined. Here, we ask:

\medskip

\begin{quote}
Can an unknown linear combination of $ N $ exponential functions with unknown rate constants be uniquely determined from its sampled values at any $ 2N $ points? That is, given
\begin{equation}\label{Eq:SOE}
    f_N(t) = \sum_{n = 1}^N c_n e^{\alpha_n t}
\end{equation}
can we determine $ c_n $ and $ \alpha_n $, $ n = 1, 2, \cdots N$ from the values of $ f(t)$ at $ t = t_1, \cdots , t_{2N} $? 

\end{quote}

The motivation behind this simple statement is quite significant. In many systems of linear differential equations the solution is a linear combination of exponentials, where the rate constants $ \alpha_i $ and the coefficients $ c_i $ (eigenvalues and components of the eigenvectors of the system) are unknown. Finding an exact solution to \eqref{Eq:SOE} provides a method to solve the inverse problem: Data as a function of time gathered at $ 2N $ data points may be used to solve for $ \alpha_i $ and $ c_i $, $ i = 1, \cdots, N $. Other applications of this problem to mixing theory in probability and statistics and elsewhere can be cited \cite{doi:10.1080}. To relate this problem to the fundamental theorem of algebra, note that under the change of variable $ x = e^t $, \eqref{Eq:SOE} takes the form
\begin{equation}\label{Eq:series}
f(x) = \sum_{n=1}^N c_n x^{\alpha_n} 
\end{equation} 
where $ \alpha_n \in \R $. Here, the exponents are no longer restricted to the positive integers. Thus, the stated problem is one of finding not only the coefficients, but the powers of a series of the form \eqref{Eq:series} from the values of $ f (x_i)  $ at $\{ x_i > 0 :  i = 1, 2, \ldots , 2N \} $. 

\medskip

This problem has received significant attention over the past 50 years by well known numerical and computational analysts (for example, see \cite{golub1973differentiation}, \cite{pereyra2006variable}, \cite{etna_vol34_pp163-169}, \cite{mca27040060}, and the many references therein). Here we propose an interesting and simple algorithm that seems to have gone unnoticed. Some key results in harmonic analysis involving Fourier series, Shannon's sampling theorem \cite{1455040} and the Poisson summation formula \cite{enwiki:1187353254} deal with closely related questions. This methodology may prove useful in identifying linear combinations of exponentials in terms of discrete set of \textit{a priori} unknown frequencies. Moments of \eqref{Eq:SOE} over varying intervals make a cameo appearance in this development.

%\medskip

%Initially, we will assume that the number of exponentials $ N $ in Eq. \eqref{Eq:SOE} is known. In a later section, we will loosen this assumption slightly to where  the problem may become under-determined/over-determined (i.e. fewer/greater than $2N$ data points are given for a linear combination of $N$ exponentials).

\section{The $N = 1$ and $ N=2 $ cases}
To illustrate the proof of the theorem, and the additional data that will be needed to solve for the unknown parameters, we will first present the cases $ N=1 $ and $ N=2 $ explicitly. These examples should at once make the idea of the proof of the general case clear. Having these models at hand, we shall state and prove the main theorem of the paper and discuss the existence, uniqueness and some variations on the problem. Extensions and applications of these results will be given elsewhere.

\medskip

For non-degeneracy of the problem, let's assume that the rate constants are nonzero and distinct. We shall deal with the special case of the zero rate constant (constant term) afterwards.

\subsection{The $N=1$ case}
Let $ 0 < t_1 < t_2 < \infty $ and in addition to knowing $ f_1 (t_1) $ and $ f_1 (t_2) $ suppose that we know the partial cumulative distribution functions $ \int_0^{t_i} f_1 (t') dt' $ for $ i = 1,2 $. Clearly these extra data are related to the average of $ f $ on $ [0,t_i] $. Instead of the traditional means of fitting an exponential to two points, note that if $ f_1 (t) = ce^{\alpha t} $ then
\begin{equation}\label{2}
    I_0 := \int_0^t f_1(t') dt' = \frac{c}{\alpha}e^{\alpha t} - \frac{c}{\alpha}.
\end{equation}
We also note that \eqref{2} is the $ 0$-th moment of $ f_1 $ on $ [0,t] $. There exists some $x_1$, $x_2$ such that
\begin{equation}
    x_1\left(\frac{c}{\alpha}e^{\alpha t} - \frac{c}{\alpha}\right) + x_2 = c e^{\alpha t}, 
\end{equation}
\begin{equation}
    x_1\int_0^t f_1(t') dt' + x_2 = c e^{\alpha t} = f_1(t).
\end{equation}
This points us to consider this equation in vector formulation.
\begin{equation}
    f_1(t) = c e^{\alpha t} = \begin{bmatrix}
        \int_0^t f_1(t') dt' & 1
    \end{bmatrix}
    \begin{bmatrix}
        x_1\\x_2
    \end{bmatrix}.
\end{equation}
Here we have shifted from needing two points to needing to know the value of the integral at two $t$ values. So with two points $ t_1 $ and $ t_2 $ we have

%For practical applications, this is trivial since one of several numerical integration techniques can be used to evaluate the integral at points. 
\begin{equation}\label{Eq:System}
    \begin{bmatrix}
        f_1(t_1)\\f_1(t_2)
    \end{bmatrix}=\begin{bmatrix}
        \int_0^{t_1} f_1(t') dt' & 1 \\
        \int_0^{t_2} f_1(t') dt' & 1
    \end{bmatrix}
    \begin{bmatrix}
        x_1\\x_2
    \end{bmatrix} .
\end{equation}
Since, by hypothesis,  $f_1(t_1)$ and $f_1(t_2)$ and the averages of $ f$ on $[0,t_i] $ are known, we have a linear system of equations. If $ t_1 < t_2 $, by positivity of $ f_1 (t) $ (if $ c > 0 $ or the reverse inequality if $ c < 0 $), $ \int_0^{t_1} f_1(t')dt' < \int_0^{t_2} f_1(t')dt' $. Hence the above $ 2 \times 2 $ matrix has a strictly negative determinant and is invertible, and the system in \eqref{Eq:System} will have a unique solution.

%The least square solution is easily found for $>2N$ points.
Recalling \eqref{2}, and substituting it for the integral, we have the system of equations
\begin{align}
    c e^{\alpha t_1} &= x_1\left(\frac{c}{\alpha}e^{\alpha t_1} - \frac{c}{\alpha}\right) + x_2 \\
    c e^{\alpha t_2} &= x_1\left(\frac{c}{\alpha}e^{\alpha t_2} - \frac{c}{\alpha}\right) + x_2,
\end{align}
which gives
\begin{align}
    x_1 &= \alpha\\
    x_2 &= c.
\end{align}

The $N = 1$ case seems like a circuitous method to go about parameter fitting, because it is. Other cases demonstrate the power of this technique.

\subsection{The $N = 2$ case}

Using the same set up as in the $N =1$ case, we see how a completely unmodified approach is ineffectual. This is due to greater number of basis functions requiring a greater number of unknowns to compensate for their integrals. Let
$ \alpha_1 \neq \alpha_2 $, 
\begin{equation}
    f_2(t) = c_1 e^{\alpha_1 t} + c_2 e^{\alpha_2 t},
\end{equation}
and assume that $ f_2 (t) \geq 0 $ for $ t \geq 0 $ (this artificial assumption will be considered in Corollary 4). As in the $N=1 $ case, assume that we also know that partial cumulative distribution functions $ I_0 = \int_0^{t_i} f_2(t') dt' $ and $ I_1 := \int_0^{t_i} \int_0^{t'} f_2(t^{\prime \prime}) dt^{\prime\prime} dt^\prime $ for $ i=1,2,3,4$. While the single integral of $ f_2 $ up to $ t_i $ is the $ 0 $-th moment, $ m_0 $, of $ f_2 $ on $ [0,t_i]$, the double integral is related to the $0$-th and $ 1$-st moment, $ m_1(t) $, of $ f_2 (t) $ on $ [0,t] $ as $  m_0 (t) t - m_1 (t) $. Then
\begin{equation}
    \int_0^t f_2(t')dt' = \frac{c_1}{\alpha_1} e^{\alpha_1 t} + \frac{c_2}{\alpha_2} e^{\alpha_2 t} - \frac{c_1}{\alpha_1} - \frac{c_2}{\alpha_2}.
\end{equation}
Since $\alpha_1 \neq \alpha_2$, we cannot draw the same equivalence we did for $N = 1$. Taking repeated integral of $f_2(t)$ we have
\begin{align}
    \int_0^t\int_0^{t'} f_2(t'')dt''dt' &= \int_0^t \left(\frac{c_1}{\alpha_1} e^{\alpha_1 t'} + \frac{c_2}{\alpha_2} e^{\alpha_2 t'} - \frac{c_1}{\alpha_1} - \frac{c_2}{\alpha_2}\right) dt'\\
    &= \frac{c_1}{\alpha_1^2} e^{\alpha_1 t} + \frac{c_2}{\alpha_2^2} e^{\alpha_2 t} - \left(\frac{c_1}{\alpha_1} + \frac{c_2}{\alpha_2}\right)t - \left(\frac{c_1}{\alpha_1^2} + \frac{c_2}{\alpha_2^2}\right).
\end{align}
With the repeated integration alone, we still cannot draw the same equivalence. However, a linear combination of the first and second integrals along with terms to compensate for constant and $t$ terms, we have that there exists some $x_1, x_2, x_3, x_4$ such that
\begin{align}\label{Eq:linear_comb}
    f_2(t) &= x_1 \int_0^t f_2(t')dt' + x_2 \int_0^t\int_0^{t'} f_2(t'')dt''dt' + x_3 t + x_4\\
    &= \begin{bmatrix}\label{Eq:System2}
        \int_0^t f_2(t')dt' & \int_0^t\int_0^{t'} f_2(t'')dt''dt' & t & 1
    \end{bmatrix}\begin{bmatrix}
        x_1\\x_2\\x_3\\x_4
    \end{bmatrix}.
\end{align}
Evaluating \eqref{Eq:System2} at $ t_1, t_2, t_3 $ and $ t_4 $, we arrive at a system of $ 4 $ linear equations in the $ 4 $ unknowns $x_1, x_2, x_3, x_4$. By positivity of $ f_2 (t) $, it is easy to verify that this $ 4 \times 4 $ matrix has a strictly negative determinant and hence is invertible, thus giving a unique solution $x_1, x_2, x_3, x_4$.

\medskip

Now we solve for $c_1, c_2, \alpha_1$ and $\alpha_2$ in terms of $ x_1, x_2, x_3 $ and $ x_4$.

\begin{multline}
c_1 e^{\alpha_1 t} + c_2 e^{\alpha_2 t} = x_1\left(\frac{c_1}{\alpha_1} e^{\alpha_1 t} + \frac{c_2}{\alpha_2} e^{\alpha_2 t} - \frac{c_1}{\alpha_1} - \frac{c_2}{\alpha_2}\right) \\
   + x_2\left(\frac{c_1}{\alpha_1^2} e^{\alpha_1 t} + \frac{c_2}{\alpha_2^2} e^{\alpha_2 t} - \left(\frac{c_1}{\alpha_1} + \frac{c_2}{\alpha_2}\right)t - \left(\frac{c_1}{\alpha_1^2} + \frac{c_2}{\alpha_2^2}\right)\right)
   + x_3 t + x_4.
\end{multline}
Terms that are constants will cancel, as will terms which are multiples of $t$.
%______________________________________
% Based on this equation,  why don't we want to just take the 2x2 system with the first 2 terms of the matrix, and $ x_1 $ and $ x_2 $.  That should be sufficient to solve for alpha_1 and alpha_2. 
%_____________________________________
\begin{equation}\label{Eq:Cancel}
c_1 e^{\alpha_1 t} + c_2 e^{\alpha_2 t} = x_1\left(\frac{c_1}{\alpha_1} e^{\alpha_1 t} + \frac{c_2}{\alpha_2} e^{\alpha_2 t}\right) 
   + x_2\left(\frac{c_1}{\alpha_1^2} e^{\alpha_1 t} + \frac{c_2}{\alpha_2^2} e^{\alpha_2 t}\right).
\end{equation}
From this we have the linear system using the exponential functions as the basis
\begin{align}
    1 &= \frac{x_1}{\alpha_1} + \frac{x_2}{\alpha_1^2}\\
    1 &= \frac{x_1}{\alpha_2} + \frac{x_2}{\alpha_2^2}.
\end{align}
It is interesting to note that due to the above cancellations, $ x_1 $ and $ x_2 $ will determine $ \alpha_1 $ and $ \alpha_2 $ uniquely. Considering one of these equations, we arrive at a quadratic equation for $ \alpha_1$. 
\begin{align}
    1 = \frac{x_1}{\alpha_1} + \frac{x_2}{\alpha_1^2}\\
    \alpha_1^2 = \alpha_1 x_1 + x_2\\
    \alpha_1^2 - \alpha_1 x_1 - x_2 = 0.
\end{align}
This quadratic equation can be solved giving us two values, one for each $\alpha $.
\begin{equation}
    \alpha = \frac{x_1 \pm \sqrt{x_1^2 + 4x_2^2}}{2}.
\end{equation}
Having $ \alpha_1 $ and $ \alpha_2 $, the coefficients $ c_1 $ and $ c_2 $ can be readily determined from a linear system of equations. 

\medskip

A simple but important point may require clarification here. We note that $ \alpha_1 $ and $ \alpha_2 $ seem to only depend on $ x_1 $ and $ x_2 $, so it is tempting to rewrite \eqref{Eq:linear_comb} as a linear of combination of the first two terms alone, and to dismiss $ x_3 $ and $ x_4 $ altogether. However, doing so will not allow the cancellation in \eqref{Eq:Cancel} to occur, and would leave that equation as indeterminate. 

\section{A fundamental theorem of algebra for \\exponentials: Existence and Uniqueness}

The above special cases have now paved the way to state and prove the general case. As mentioned earlier, for nondegeneracy assume $ \alpha_i \neq \alpha_j $ for $ i \neq j $; otherwise, \eqref{Eq:SOE} reduces to a linear combination of $ N-1 $ or fewer terms. Initially, we shall assume that the function $ f_N (t) $ in \eqref{Eq:SOE} is nonnegative for all $ t \geq 0 $ and $ \alpha_n \neq 0 $ for all $ n $. We shall see later that these assumptions can be dismissed via simple translation. 

\medskip

Let $ f_N(t) $ be given by \eqref{Eq:SOE} and define
\begin{equation}\label{Eq:II}
    I_{(k)}(t) := \int_0^{t}\cdots\int_0^{t^{(k)}}  f_N(t^{(k+1)}) dt^{(k+1)} \cdots dt^{(1)}, \ k=0,1,\cdots 
\end{equation}
where $ t^{(k)} $ is used to denote the $ k $-th variable of integration, $ t^{(0)} = t $, and there are $ k+1 $ integrals. As it noted in the $ N=2 $ case, we can calculate the succeeding parameters $ I_{(k)} (t) $ in terms of the first $ k$ moments of $ f_N (t) $ on $[0,t]$. Direct calculation shows that 
$$
\begin{array}{ll}
 I_{(1)} (t) = &=  m_0 (t) t - m_1 (t) \\
 2! \ I_{(2)}(t)  &= m_0(t) t^2 - 2 m_1(t) t  + m_2(t) \\
 3! \ I_{(3)} (t) &= m_0 (t) t^3  - 3 m_1 (t) t^2  + 3 m_2 (t) t  - m_3 (t)(t) .
 \end{array}
$$
In general, the relationship between the moments on $ [0,t] $ and the iterated integrals $ I_{(k)} (t) $ is given by

\noindent{\bf Proposition 1.}
{\it Let $ I_{(k)} (t) $ be described by \eqref{Eq:II} and $ m_k (t) $ be the $ k $-th moment of $ f_N (t) $ in \eqref{Eq:SOE} on $ [0,t] $. Then
\begin{equation}\label{binom}
I_{(k)} (t) = \frac{1}{k!} \sum_{j=0}^k (-1)^j \binom{k}{j} \ m_j (t) \ t^{k-j} .
\end{equation}
}

\noindent{\bf Proof.}
By linearity, it is sufficient to show \eqref{binom} for each term in $ f_N (t)$, namely for $ e^{\alpha t} $. Note that 
\begin{equation}
    I_{(k)} (t) = \int_0^t I_{(k-1)} (t') dt'
\end{equation}
and by integation by parts, 
\begin{equation}
    m_k (t) = \int_0^t {t^\prime}^k e^{\alpha t^\prime} dt^\prime  =  \frac{t^k e^{\alpha t}}{\alpha} - \frac{k}{\alpha} m_{k-1} (t).
\end{equation}
Then a straight forward induction argument shows that 
\begin{equation}
k! I_{(k)} (t) - (-1)^{k} m_k(t) = \sum_{j=0}^{k-1} (-1)^j \binom{k}{j} \ m_j (t) \ t^{k-j}.
\end{equation}
\hfill$\Box$
\medskip

\noindent{\bf Theorem 2.}
   {\it  Let $ f_N(t) = \sum_{n = 1}^N c_n e^{\alpha_n t} $ and assume that $ f_N (t) \geq 0 $ and $ \alpha_n \neq 0 $ for all $ n $. Given $\{ f_N(t_i): i = 1, \cdots 2N\}  $, and $ \{I_{(k)} (t_i):  k = 0, 1, \cdots, N, i = 1, \cdots, 2N\} $, then $ c_n $ and $ \alpha_n $ are uniquely determined for $ n \in \{1,2, \cdots, N\} $.}

\medskip

\noindent{\bf Proof.}
We generalize the repeated integration of $f_N(t)$ that was demonstrated in the case  $N = 2$, extending it $ N $ times. Order $ t $'s such that $ 0 \leq t_1 < t_2 < \cdots < t_{2N} < \infty $. Using $ I_k (t) $ defined in \eqref{Eq:II}, inserting $ f_N (t) $, and integrating $ k $ times, we obtain
\begin{equation}
    I_{(k-1)}(t) = \sum_{n=1}^N \left(\frac{c_n}{\alpha_n^k} e^{\alpha_n t} - \sum_{j=0}^{k - 1}\frac{c_n}{\alpha_n^{k-j}j!}t^{j}\right).
\end{equation}
We write $f_N(t)$ as a linear combination of $t^n $ and $I_{(k-1)}(t)$\ 's  with coefficients $x_1, x_2,\ldots, x_{2N}$, to construct the prototypical row of a $ 2N \times 2N $ linear system of equations
\begin{equation}\label{Eq:System_N}
    f_N(t) = \begin{bmatrix}
        I_{(0)}(t) & I_{(2)}(t) & \cdots & I_{(N-1)}(t) & 1 & t & t^2 & \cdots & t^{N-1}
    \end{bmatrix}\begin{bmatrix}
        x_1\\ \vdots \\ x_{2N}
    \end{bmatrix}
\end{equation}
or, more explicitly,
\begin{equation}
    f_N(t) = x_1 I_{(1)}(t)+x_2I_{(2)}(t) + \cdots + x_NI_{(N)} (t)+ x_{N+1} + x_{N+2}t + \cdots + x_{2N}t^{N-1}.
\end{equation}
This system of equations is obtained by evaluating each of the terms at $ t_1, \cdots, t_{2N} $. By the positivity of $ f_{N} (t) $, the $ 2N \times 2N $ matrix in \eqref{Eq:System_N} has a nonzero determinant (proved by induction), and hence the system of equations has a unique solution.

\medskip

Distributing the $x_k$ coefficients within each $I_{(k)}(t)$, we have that
\begin{align}
    x_k I_{(k-1)}(t) &= x_k\sum_{n=1}^N \left(\frac{c_n}{\alpha_n^k} e^{\alpha_n t} - \sum_{j=0}^{k - 1}\frac{c_n}{\alpha_n^{k-j}j!}t^{j}\right)\\
    &= x_k\sum_{n=1}^N \frac{c_n}{\alpha_n^k} e^{\alpha_n t} - x_k\sum_{n=1}^N\left(\sum_{j=0}^{k - 1}\frac{c_n}{\alpha_n^{k-j}j!}t^{j}\right)\\
    &= x_k\sum_{n=1}^N \frac{c_n}{\alpha_n^k} e^{\alpha_n t} - \sum_{j=0}^{k - 1}\left(x_k\sum_{n=1}^N\frac{c_n}{\alpha_n^{k-j}j!}t^{j}\right).
\end{align}
Therefore the sum of $t^j$ terms can be expressed as
\begin{equation}
    x_{N+1+j}t^j-\sum_{k = 1}^{N}x_k\sum_{n=1}^N\frac{c_n}{\alpha_n^{k-j}j!}t^{j} = 0
\end{equation}
with equality to zero since $f_N(t)$ has no $t^n$ terms.

Further, the sum of $c_m e^{\alpha_m t}$ terms can be expressed as
\begin{equation}
    \sum_{k=1}^N x_k \frac{c_m}{\alpha_m^k}e^{\alpha_m t} = c_m e^{\alpha_m t}
\end{equation}
\begin{equation}
    \sum_{k=1}^N x_k \frac{1}{\alpha_m^k} = 1
\end{equation}
%\begin{equation}
%    c_m e^{\alpha_m t}\sum_{k=1}^N x_k \frac{1}%{\alpha_m^k} = c_m e^{\alpha_m t}
%\end{equation}
%\begin{equation}
%    \sum_{k=1}^N x_k \frac{1}{\alpha_m^k} = 1
%\end{equation}
%\begin{equation}
%    0 = 1 - \sum_{k=1}^N x_k \frac{1}%{\alpha_m^k}
%\end{equation}
%\begin{equation}
%    0 = \alpha_m^N - \sum_{k=1}^N x_k %\frac{\alpha_m^N}{\alpha_m^k}
%\end{equation}
\begin{equation}
    \alpha_m^N - \sum_{k=1}^N \alpha_m^{N-k}x_k = 0
\end{equation}
\begin{equation}\label{Frobenius}
    \alpha^N - \sum_{k=1}^N \alpha^{N-k}x_k = 0
\end{equation}
which fully determines $\{\alpha_n\}$. It is noteworthy that the $ \alpha_i $ are the different roots of  \eqref{Frobenius}, a monic Fr\"{o}benius polynomial, and hence they are the eigenvalues of the corresponding companion matrix. Since we assumed (for non-degeneracy of the original equation) that $ \alpha_i \neq \alpha_j $, the Fr\"{o}benius polynomial will have distinct roots, and the companion matrix is diagonalizable with a Vandermonde matrix whose $ i$th row is $ [1 \ \alpha_i \ \alpha_i^2 \ \cdots \ \alpha_i^N] $.\

\hfill$\Box$

\medskip

What is occurring here is a change of basis in the expression of the original equation. Changing from the basis consisting of $I_{(k)}(t)$ and powers of $t$ to the basis consisting of exponential functions and powers of $t$. The existence of the vector $\mathbf{x}$ from the first basis and following the change of basis, this vector gives us the values of the vector $ \mathbf{\alpha}$. 

\medskip

\noindent{\bf Proposition 3.} {\it If $ r \geq N $ and $ t_{2N+1}, \cdots , t_{2r} $ are additional points where $ f_N (t_i) $ are specified, then $ c_i $ and $ \alpha_i $ will be the same for $ i = 1, \cdots N$ and $ c_{N+1}, \cdots , c_r $ will be zero.}

\medskip

\noindent{Proof.} Without loss of generality, order the specified $ t $'s such that $ \{t_1, \cdots , t_{2N} \}$ appear first in the list and $ t_{2N+1}, \cdots , t_{2r} $ are the additional points. Theorem 1 already demonstrates that given $ 2N $ $ t $'s (counting multiplicity), the $ c_i $ and $ \alpha_i $ are determined uniquely. 
Hence, suppose that \eqref{Eq:SOE} is determined from $\{t_1, \cdots, t_{2N} \}$ and
\begin{equation}\label{Eq:SOE-plus}
    f_N(t) = \sum_{n = 1}^{r} c_n e^{\alpha_n t}
\end{equation}
Subracting \eqref{Eq:SOE} from \eqref{Eq:SOE-plus}, we get
\begin{equation}\label{Eq:Indep}
0 = \sum_{n = N+1}^{r} c_n e^{\alpha_n t}.    
\end{equation}
Calculating the Wronskian of the set of exponentials $ \{ e^{\alpha_{N+1} t}, \cdots, e^{\alpha_{r} t} \} $, since $ \alpha_i $ are assumed to be distinct, the resulting Vandermonde matrix is strictly positive definite, hence \eqref{Eq:Indep} forces $ c_i = 0 $ for $ i = N+1, \cdots, r $.
\hfill$\Box$

\medskip

The two (artifical) assumptions made in Theorem 2 require brief discussion.  We assumed that $ \alpha_i \neq 0 $ for any $ i $ and $ f_N(t) \geq 0 $ for all $ t \geq 0 $.  Clearly, if $ \alpha_i = 0 $ for some $ i $, say $ i=1 $, then that term is constant, and  \eqref{Eq:SOE} reduces to
\begin{equation}
f_N(t) = c_1  + \sum_{i=2}^N c_i e^{\alpha_i t }.
\end{equation}
Hence there will be only $ (N-1) $, $\alpha$'s, and $ N $ $c_i $'s to be determined. Once the $\alpha $'s are determined, the system is linear in $ c_i $, and those can be readily determined. Here $ (2N - 1) $, $ t_i $'s will be sufficient. 

%________________________________
%\textcolor{red}{Please check that nothing goes wrong here. The N=2 case is sufficient for this check.}

%\textcolor{blue}{I have N=2 worked out on paper where $\alpha_1 = 0$.}
%________________________________

\medskip

The assumption of $ f_N(t) \geq 0 $ is somewhat more intricate. This assumption originally appeared in \eqref{Eq:System2} to guarantee the invertibility of the $ 2 \times 2 $ matrix in that equation.  This nonnegativity also readily allows the invertibility of the $ 2N \times 2N $ matrix in  \eqref{Eq:System_N}. Clearly, this is a suficiency result. For a general $ f $ (not necessarily nonnegative), this assumption may be replaced by choosing $ \{t_i: i = 1, \cdots, 2N \} $ so that the $ 2N \times 2N $ matrix in \eqref{Eq:System_N} is invertible. Alternatively, this can be treated via dimension extension. Introducing $ \alpha_{n+1} = 0 $ and changing \eqref{Eq:SOE} to 
\begin{equation}\label{Eq:dimext}
f_{N+1} (t) = c_{N+1} + \sum_{n=1}^N c_n e^{\alpha_n t} ,
\end{equation}
if $ c_{N+1} $ is chosen sufficiently large (e.g. $ c_{n+1} \geq \min\{f_N(t); t \in [0,T] \} $), $ f_{N+1} (t) $ will be $ \geq 0 $ for all $ t \in [0,T] $. Then with $ 2N+1 $ data points, \eqref{Eq:dimext} is solveable for $ \alpha $'s and $ c $'s.  We summarize this observation in the following and provide the details of a proof.

\medskip

\noindent{\bf Corollary 4.} 
{\it  Let $ f_N(t) = \sum_{n = 1}^N c_n e^{\alpha_n t} $. Given $\{ f_N(t_i): i = 1, \cdots 2N \}  $, and $ \{I_{(k)} (t_i):  k = 1, 2, \cdots, N, i = 1, \cdots, 2N+1 \} $, then $ c_i $ and $ \alpha_i $ for $ i \in \{1,2, \cdots, N\} $, are uniquely determined.}

\medskip

\noindent{\bf Proof.}
As mentioned above, consider the closely related problem
\begin{equation}
f_{N+1} (t) = c_{N+1} + \sum_{n=1}^N c_n e^{\alpha_n t} = c_{N+1} + f_N (t) .
\end{equation}
where $ c_{N+1 } \geq \min\{f_N (t) : t \in [0,T]\} $. Let $ \tilde{I_k} (t) = I_k (t) + \frac{c_{N+1}t^k}{k!} $ for all $ k \in \{1, \cdots. N\}$. Then $ f_{N+1}(t)  \geq 0 $ for $ t \in [0,T] $. Using Theorem 2, we can determine $ \{ \alpha_n : n = 1, \cdots, N\} $. This will give a linear system in $ \{c_n: n=1, \cdots N+1\} $.  Hence the problem is uniquely solvable.

\hfill$\Box$

\section{Concluding Remarks}

We conclude with a few immediate observations of the results in this paper. Similarly to the classical fundamental theorem of algebra, the choice of $ \{t_i: i = 1, \cdots, 2N \} $ is completely arbitrary in determining $ c_i $ and $ \alpha_i $ in \eqref{Eq:SOE}. No assumption on uniform distribution of these points is required. This observation has major consequences in the sampling problem and in the inversion problem in the differential equation and applied settings. The appearance of the cumulative integrals or the moments of $ f_N(t) $ up to $ [0,t_i] $, $ i = 1, \cdots , 2N $, introduces sampling and statistical considerations in applied problems, where data may be quantized. Expressing these culmulative integrals as polynomials in the moments as it is done in Proposition 1 opens a new connection to be exploited much further. The vast literature on moments (see \cite{akhiezer2020classical}, \cite{Dragotti} and the references therein) should 
be useful.

\medskip

A natural extension of this technique is fitting data to sums of exponentials where we do not know the exact number of exponential basis functions. The over-determined problem has already been considered in Proposition 2, while the 
under-determined problem will require optimization over the set $ c_i $ and $ \alpha_i $. Many interesting questions will arise in this context, whose treatment will require significant space and numerical analysis. The interesting works \cite{Dragotti} and \cite{etna_vol34_pp163-169} have already established least-square approaches to this problem. The present work should provide nice extensions to their methods.

\providecommand{\bysame}{\leavevmode\hbox to3em{\hrulefill}\thinspace}
\providecommand{\MR}{\relax\ifhmode\unskip\space\fi MR }
\providecommand{\MRhref}[2]{%
  \href{http://www.ams.org/mathscinet-getitem?mr=#1}{#2}
}
\providecommand{\href}[2]{#2}

% Bibliography %------------------------------%
%\bibliographystyle{amsplain}
%{\small\bibliography{references.bbl}}

\end{document}